\documentclass[11pt]{amsart}

\usepackage{amsthm}
\usepackage{amssymb}
\usepackage{amsmath}

\newcommand{\labell}[1] {\label{#1}}

\numberwithin{equation}{section}

\newtheorem {Theorem}   {Theorem} 
\numberwithin{Theorem}{section}

\theoremstyle{definition}

\theoremstyle{remark}
\newtheorem{Remark}[Theorem]{Remark}
\newtheorem{Example}[Theorem]{Example}

\newtheorem {Corollary}[Theorem]{Corollary}  

\expandafter\chardef\csname pre amssym.def at\endcsname=\the\catcode`\@ 
\catcode`\@=11 
\def\undefine#1{\let#1\undefined} 
\def\newsymbol#1#2#3#4#5{\let\next@\relax 
 \ifnum#2=\@ne\let\next@\msafam@\else 
 \ifnum#2=\tw@\let\next@\msbfam@\fi\fi 
 \mathchardef#1="#3\next@#4#5}
\def\mathhexbox@#1#2#3{\relax 
 \ifmmode\mathpalette{}{\m@th\mathchar"#1#2#3}% 
 \else\leavevmode\hbox{$\m@th\mathchar"#1#2#3$}\fi} 
\def\hexnumber@#1{\ifcase#1 0\or 1\or 2\or 3\or 4\or 5\or 6\or 7\or 8\or 
 9\or A\or B\or C\or D\or E\or F\fi} 
 
\font\teneufm=eufm10 
\font\seveneufm=eufm7
\font\fiveeufm=eufm5 
\newfam\eufmfam
\textfont\eufmfam=\teneufm 
\scriptfont\eufmfam=\seveneufm 
\scriptscriptfont\eufmfam=\fiveeufm 
 
\catcode`\@=\csname pre amssym.def at\endcsname 

\def	\eps	{\epsilon}

\newcommand{\CE}{{\mathcal E}}

\newcommand{\const}{{\mathit const}}
\newcommand{\divv}{{\mathit div}\,}

\newcommand{\hook}{\hookrightarrow}

\def	\C {{\mathbb C}}
\def	\reals	{{\mathbb R}}
\def	\R	{{\mathbb R}}

\def	\rationals	{{\mathbb Q}}
\def	\TT	{{\mathbb T}}
\def	\T	{{\mathbb T}}
\def	\CP	{{\mathbb C}{\mathbb P}}

\def	\p	{\partial}

\def	\ssminus 	{\smallsetminus}

\begin{document}

%%%%%%%%%%%%%%%%%%%%%%%%%%%%%%
%   TEXT FORMATTING

\setlength{\smallskipamount}{6pt}
\setlength{\medskipamount}{10pt}
\setlength{\bigskipamount}{16pt}

%%%%%%%%%%%%%%%%%%%%%%%%%%

%%%%%%%%%%%%%%%%%%%%%%%%%%

%%%%%%%%%%%           BEGINNING OF  TEXT

%%%%%%%%%%%%%%%%%%%%%%%%%%

\title[Hamiltonian Dynamical Systems Without Periodic Orbits]{Hamiltonian
Dynamical Systems Without Periodic Orbits}

\author[Viktor Ginzburg]{Viktor L. Ginzburg}

\address{Department of Mathematics, UC Santa Cruz, 
Santa Cruz, CA 95064}
\email{ginzburg@math.ucsc.edu}

\date{October, 1998}

\thanks{The work is partially supported by the NSF and by the faculty
research funds of the University of California, Santa Cruz.}

\bigskip

\begin{abstract}
The present paper is a review of counterexamples to the 
``Hamiltonian Seifert conjecture'' or, more generally,
of examples of Hamiltonian systems having no periodic orbits on a 
compact energy level. We begin with the discussion of 
the ``classical'' and volume--preserving Seifert conjectures.
Then a construction of counterexamples to the Hamiltonian
Seifert conjecture in dimensions greater than or equal to six is
outlined, mainly following the method introduced by the author
of the paper. The Hamiltonian version of Seifert's theorem is also
stated.  A list of all known at this moment examples and 
constructions of Hamiltonian flows without periodic flows concludes
the review.
\end{abstract}

\maketitle

\section{Introduction} 

The existence problem for periodic orbits of vector fields or diffeomorphisms
occupies one of the central places in the theory of dynamical systems and 
adjacent areas such as mechanics and symplectic geometry.
This question is often well motivated for a particular system,
but in many cases the answer gives very little information about 
the dynamics of the system in general. However, the methods developed 
to solve the existence problem may have an impact extending well beyond 
the scope of the original problem.

The search for dynamical systems without periodic orbits has been inspired
by a few questions. One of them is to determine the 
limits of the existence theorems. For instance,
Seifert's theorem on periodic orbits of vector fields on $S^3$
led to the famous Seifert conjecture recently disproved by K. Kuperberg,
\cite{kuk}. However, in addition to this, the systems found 
as a result of this search
sometimes exhibit a new type of dynamics and extend our understanding 
of the qualitative behavior that can occur for a given class 
of flows.

In this review we focus on examples of Hamiltonian systems without
periodic orbits on a compact energy level in the context of the
Seifert conjecture. The paper is organized
as follows. 

In Section \ref{sec:seifert} we recall Seifert's theorem
and the Seifert conjecture. We discuss the history of counterexamples
to the Seifert conjecture from Wilson's theorem (with its proof)
to the ultimate solution due to K. Kuperberg. We also briefly touch upon
related results on volume--preserving flows.

Section \ref{sec:Ham} is devoted to Hamiltonian vector
fields. We outline the constructions of counterexamples to the Hamiltonian
Seifert conjecture in dimensions greater than or equal to six,
mainly following the method introduced by the author of the 
present review. We also state the Hamiltonian version of Seifert's 
theorem.

Finally, a list of all known at this moment 
constructions of Hamiltonian flows without periodic flows is
presented  in Section~\ref{sec:list}.

\subsection*{Acknowledgments.} The author is deeply grateful to
Ana Cannas da Silva, Carlos Gutierrez, Ely Kerman, Greg and Krystyna
Kuperberg, Debra Lewis,
Richard Montgomery, Marina Ratner, Andrey Reznikov, Claude Viterbo,
and the referee for their advice, remarks,
and useful discussions. 
He would also like to  thank the Universit\'{e} Paris-Sud 
for its hospitality during the period when the work on this manuscript 
was started.

\section{The Seifert Conjecture}
\labell{sec:seifert}

\subsection{Seifert's Theorem}
The history of examples of dynamical systems without periodic orbits
as we understand it today begins with a result of Seifert that
a $C^1$-smooth vector field on $S^3$ which is $C^0$-close to the Hopf
field has at least one periodic orbit, \cite{Seifert-1950}. Later this
theorem was generalized by Fuller, \cite{Fuller}, as follows.

\begin{Theorem} 
\labell{Theorem:fuller}
Let $E\to B$ be a principal circle bundle over a compact manifold 
$B$ with $\chi(B)\neq 0$. Let also $X$ be a $C^1$-smooth vector field 
$C^0$-close to the field $X_0$ generating the $S^1$-action on $E$. Then 
$X$ has at least one periodic orbit.
\end{Theorem}

Today we know two approaches to proving existence theorems for
periodic orbits such as Theorem \ref{Theorem:fuller}, neither of which 
has been trivialized and made into a part of mathematical pop culture.
The first approach, not counting the original Seifert proof, relies on
the notion of the Fuller index, \cite{Fuller}, an analogue of the Euler 
characteristic for periodic orbits. The second one, analytical, is 
due to Moser \cite{Moser-1976}. Moser's method uses a version of an 
infinite-dimensional inverse function theorem or, more precisely, its 
proof. In fact, as has been recently noticed by Kerman, \cite{ely}, 
under the additional assumption that $X$ is $C^1$-close to $X_0$ 
(but not only $C^0$-close), Moser's argument can be significantly 
simplified by applying the inverse function theorem in Banach spaces.
This stronger closeness hypothesis holds in virtually 
all applications of the theorem known to the author. Note also that 
in results such as Theorem~\ref{Theorem:fuller}, bounding
$X-X_0$ with respect to a higher order norm
often considerably simplifies the proof.

A representative corollary of Theorem~\ref{Theorem:fuller} (for which
the $C^1$-closeness assumption is sufficient) is as follows:
\emph{
Let $f\colon \reals^{2n}\to \reals$ be a smooth function
having a non-degenerate minimum at the origin and, say, $f(0)=0$.
Assume that all eigenvalues of $d^2 f(0)$ are equal.
Then for every sufficiently small $\epsilon>0$ the level $\{f=\eps\}$
carries a periodic orbit of the Hamiltonian flow of $f$.}
In effect, the assumption that all eigenvalues of $d^2 f(0)$ are equal is
purely technical and immaterial. Moreover, there are at least $n$ 
periodic orbits of the Hamiltonian flow on $\{f=\eps\}$. 
(See \cite{weinstein-1973}, \cite{Moser-1976}, \cite{Bottkol}.) 
For example, when all eigenvalues are distinct the existence of
$n$ periodic orbits readily follows from the inverse function theorem. 

\begin{Remark} It is interesting to point out that Seifert's theorem
on $(4n+1)$-dimensional spheres (i.~e., Theorem~\ref{Theorem:fuller}
for the Hopf fibration $S^{4n+1}\to \CP^{2n}$) can be proved by the
standard algebraic topological methods. Namely, let 
$D_x\subset S^{4n+1}$ with $x\in S^{4n+1}$ be a small embedded 
$4n$-dimensional disc transversal to
the fibers and centered at $x$. We can choose $D_x$ to depend 
smoothly on $x$. Let $P(x)$ be the first intersection with $D_x$ of the 
integral 
curve of $X$ through $x$. Clearly $P(x)=x$ if and only if the integral
curve closes up after one revolution along the fiber. There exists
a unique vector $v_x$ tangent to $D_x$ at $x$ such that $P(x)$ lies
on the geodesic in $D_x$ beginning at $x$ in the direction $v_x$
and the distance from $x$ to $P(x)$ is $\parallel v_x \parallel$.
Hence, $v_x=0$ if and only if $P(x)=x$. Thus it suffices to prove
that the vector field $v$ vanishes at least at one point on $S^{4n+1}$.
Note also that $v$ is normal to $X_0$ with respect to a suitably
chosen metric and so $v$ and $X_0$ would be linearly independent if
$v$ did not vanish.

The sphere $S^{4n+1}$ does not admit two linearly independent vector 
fields. (This is a very particular case of Adams' theorem that can be
proved by using the Steenrod squares; see \cite{steenrod}.) 
Therefore, on $S^{4n+1}$ the field $v$ vanishes somewhere and $X$ 
has a periodic orbit.
\end{Remark}

\subsection{The Seifert Conjecture}
\labell{subsec:seif}
In the same paper, \cite{Seifert-1950}, where he proved his theorem
discussed above,
Seifert asked the question whether every non-singular vector field on $S^3$
has a periodic orbit or not. The hypothetically affirmative answer to this question
has become known as the Seifert conjecture. The three--dimensional
sphere plays, of course, a purely historical role in the 
conjecture and a similar question can be asked for other manifolds
and also for more restricted classes of vector fields (e.~g.,
real--analytic, divergence--free, or Hamiltonian) or for non-vanishing
vector fields in a fixed homotopy class (see Remark \ref{rmk:homotop}).

The first counterexample to the generalized Seifert conjecture is due 
to Wilson, \cite{wilson}, who showed that the smooth Seifert conjecture 
fails in dimensions greater than three by proving the following

\begin{Theorem}
\labell{thm:wilson}
Let $M$ be a compact connected manifold with $\chi(M)=0$ and
$m=\dim M\geq 4$. Then there exists a smooth vector field on $M$
without periodic orbits and singular points.
\end{Theorem}

\begin{proof}%[Outline of the proof.]
Let $v_0$ be a vector field on $M$ without zeros.
The idea is to modify $v_0$ so as to eliminate periodic orbits.
First assume that $v_0$ has a finite number of periodic orbits. 
(Note that this
is not a generic property.) Then each of the periodic orbits can be
eliminated by the following procedure which is
common to many constructions of vector fields without closed
integral curves. 

A plug is a manifold $P=B\times I$, where $I=[-1,1]$ and
$B$ is a compact manifold of dimension $m -1$ with boundary,
and a non-vanishing vector field $w$ on $P$ with the following 
properties\footnote{This definition, more restrictive than that given in
\cite{kug,kugk}, is only one of several existing definitions 
of plugs. See also \cite{kuk}, \cite{gi:seifert},
and references therein. However,
the differences between these definitions seem to be of a technical 
rather than conceptual nature.}:
\begin{enumerate}

\item 
%\labell{plug-boundary}
\emph{The boundary condition}:
$w=\p/\p t$ near $\p P$, where $t$ is the coordinate on $I$.
As a consequence, an integral curve of $w$ can only leave the plug 
through $B\times \{1\}$.

\item 
\emph{Existence of trapped trajectories}:
There is a trajectory of $w$ beginning on $B\times \{-1\}$ that
never exits the plug. Such a trajectory is said to be trapped in $P$.

\item 
\emph{Aperiodicity}:
$w$ has no periodic orbits in $P$. In other words, the ``flow'' 
of $w$ is aperiodic.

\item 
\emph{Matched ends or the entrance--exit condition}:
If two points $(x, -1)$, the ``entrance'', and $(y,1)$, 
the ``exit'', are on the same integral curve of $w$, then
$x=y$. Hence, every
trajectory of $w$ which enters and exists $P$ has its exit point right
above the entrance point (with $I$ being regarded as the vertical 
direction).

\item 
\emph{The embedding condition}:
There exists an embedding $i\colon P\hookrightarrow \reals^m$ such that 
$i_*(w)=\p/\p x_m$ near $\p P$. In other words, near the boundary
of $P$, the embedding $i$ preserves the vertical direction. 
\end{enumerate}

Assuming that the plug has been constructed, we can use it to modify $v_0$.
Namely, consider a small flow box near a point on a periodic orbit 
of $v_0$. We choose coordinates on the box so that $v_0=\p/\p x_m$.
Using $i$, we embed $P$ into the box so that the trapped trajectory
matches the periodic orbit of $v_0$. Let $v$ be the vector field
obtained by replacing $v_0$ by $w$ inside of $i(P)$. By 
the first property of the plug (the boundary condition), $v$ is 
smooth. By the third condition (aperiodicity) and
the fourth (matched ends), $v$ has one periodic orbit less than $v_0$.
By applying this procedure to every periodic orbit of $v_0$, we 
eliminate all of them.

Let us now turn to the construction of a plug. First observe
that when a plug $P'$ satisfying all of the above conditions but
the entrance--exit condition is constructed, it is easy to 
find a genuine plug $P$. Namely, $P$ is the union of two copies
of $P'$ with the second copy put upside-down on the top of the
first copy. Hence, from now on we may forget about the entrance--exit
condition. This trick, with minor modifications, is present in
all constructions of plugs.

The plug from \cite{wilson}, the so-called \emph{Wilson's plug}, is
the cylinder over the unit ball $D^{m-1}\subset \reals^{m-1}$.
Thus $B=D^{m-1}$ and $t=x_m$, and the plug is automatically embedded into 
$\reals^m$. To define $w$, fix an embedding
$\TT^2\hookrightarrow D^{m-1}\times 0$ whose normal bundle is trivial. 
(Hence the requirement that $m\geq 4$.) Then a tubular neighborhood 
$U$ of $\TT^2$ in $P$ is diffeomorphic to
$\TT^2\times D^{m-2}(\eps)$, where $D^{m-2}(\eps)$ is the
$(m-2)$-dimensional ball of a small radius $\eps>0$. Let $w_1$ be an 
irrational vector field on $\TT^2$ and let $f$ be a bump function on
$D^{m-2}(\eps)$ equal to 1 at the center of the ball and vanishing near
$\p D^{m-2}(\eps)$. Since $\TT^2$ is embedded into $D^{m-1}\times 0$, 
we may
assume that each $y\times D^{m-2}(\eps)$, $y\in\TT^2$, is parallel to
the vertical direction, i.~e., $\p/\p x_m$ is tangent to the fibers
$y\times D^{m-2}(\eps)$ of the tubular neighborhood $U$. Now we set 
\begin{equation}
\labell{eq:plug}
w=fw_1+(1-f)\frac{\p }{\p x_m}
\end{equation}
on $U$ and $w=\p/\p x_m$ on $P\ssminus U$. It is easy to check that
$(P,w)$ satisfies the plug conditions (except the entrance--exit condition).
In particular, all trapped trajectories are asymptotic to $\TT^2$. 
This completes the construction of the plug. In what
follows we will refer to $\T^2$ with an irrational flow embedded into 
the plug as to the \emph{core} of the plug. 

\begin{Remark} 
Wilson's plug defined in this section is clearly only $C^\infty$-smooth.
However, its construction can be modified to apply in the real 
analytic category; see, e.~g.,
\cite{ghys,kugk} and references therein.
\end{Remark} 

When the periodic orbits of the flow of $v_0$ are not isolated, the
plug $P$ needs to be slightly altered. Namely, one can
construct $P$ so that the beginnings in $B\times\{-1\}$
of trapped trajectories form a set with non-empty interior. Then a finite
number of plugs are inserted into $M$ so as to interrupt every 
trajectory of $v_0$ (regardless of whether it is closed or not). As before,
no new periodic orbit is created, but the original ones are eliminated.
\end{proof}

\begin{Remark}
In the plug used in the case where there are infinitely many
periodic orbits, we clearly have $\divv w\neq 0$.
This amounts to the fact that in the volume--preserving version of
Wilson's argument, one has to start with a vector field having only
a finite number of periodic orbits.
\end{Remark}

The next step after Wilson's result was the construction, due to 
Schweitzer \cite{schweitzer}, of a $C^1$-smooth non-vanishing 
vector field on $S^3$ (or on any compact three-dimensional manifold) 
without periodic
orbits. The proof again is by inserting plugs. Schweitzer's plug
uses as the core the closure of a trajectory of the Denjoy flow on 
$\TT^2$ instead of an irrational flow. (Hence, only $C^1$-smoothness.) 
Since 
the trajectory is neither dense nor closed, one can take as $B$
the torus $\TT^2$ with 
a small open disc deleted in order to have the embedding condition
satisfied. The vector field on $P$ is then given by \eqref{eq:plug}
where $w_1$ is the Schweitzer field on $B$ and $f$ is an appropriately
chosen cut-off function.
Schweitzer's construction was later improved by Harrison, 
\cite{harrison}, to obtain a $C^{2+\eps}$-smooth vector field. 

Finally, a major breakthrough came when K. Kuperberg, \cite{kuk},
constructed a $C^\infty$-smooth (and even real--analytic)
three-dimensional plug rendering thus a non-singular real--analytic
flow without periodic orbits on every compact three-manifold. 
(See also \cite{ghys} and \cite{kugk}.)
Kuperberg's plug is entirely different from those described above.
One begins with a ``plug'' $P'$ satisfying all of the conditions of the
plug but aperiodicity -- there are exactly two periodic trajectories
inside of $P'$. Then one builds $P$ by inserting parts of $P'$
into $P'$ again (self-insertion) so that the vector field on the 
inserted parts matches the original vector field. The self-insertion is
performed so as to guarantee that the resulting flow on $P$ is 
free of periodic orbits. 

\begin{Remark}
\labell{rmk:proliferation}
The common feature of the constructions described above (with
the exception of Kuperberg's plug) is that a non-singular vector field 
without periodic orbits is used as the core of the plug in order to
trap a trajectory. For example, in Wilson's plug the core flow is an 
irrational flow on the torus and in Schweitzer's plug the core is the 
Denjoy flow. Hence, one non-singular vector field without periodic orbits
gives rise to a multitude of such vector fields with various
higher--dimensional phase spaces (proliferation of aperiodic flows). 
This idea is also applied to find counterexamples to the Seifert 
conjecture in other categories. 
\end{Remark}

\begin{Remark}
\labell{rmk:homotop}
The vector field $v$ constructed in the proof of Theorem \ref{thm:wilson}
is homotopic to the original vector field $v_0$ in the class of
non-singular vector fields. This follows from that the field $w$
on Wilson's plug is homotopic to the vertical vector field. The same
holds for many other constructions of plugs including
Kuperberg's plug. Hence, the homotopy type of a non-singular vector
field is preserved while the vector field is altered by inserting
the plugs to eliminate the periodic orbits.
\end{Remark}

\subsection{Volume--Preserving Flows.}
\labell{subsec:vol}
The Seifert conjecture for this class of flows
seems to be rather similar to the Seifert conjecture in the
smooth and real analytic categories. To be more precise, as was pointed 
out by A. Katok, 
\cite{katok}, \emph{a divergence--free smooth 
non-vanishing vector field $v_0$ on a compact manifold of dimension
$m\geq 4$ can be changed into one without periodic orbits,
provided $v_0$ has only a finite number of periodic orbits.}
In fact, the field $w$ given by \eqref{eq:plug} on Wilson's plug 
(for isolated periodic orbits) can be chosen divergence--free.
This yields a smooth volume--preserving flow on $S^{2n+1}$, $2n+1\geq 5$,
without periodic orbits. 

Much less is known in dimension three. Let us state two important
results due to G. Kuperberg, \cite{kug}.

The first one is that \emph{every compact three--manifold 
$M$ possesses a volume--preserving $C^\infty$-smooth flow
with a finite number of periodic orbits and no fixed points}, 
\cite{kug}. This follows from the 
fact that $M$ can be obtained from $\T^3$ by a series of Dehn 
surgeries (provided that $M$ is orientable). Let us equip $\T^3$ with
an irrational flow. A Dehn surgery can be interpreted as the insertion
of a version of a smooth volume--preserving plug $P$ into $\T^3$. 
These plugs differ from those described in Section \ref{subsec:seif} 
in some essential ways. 
The plugs $P$ are not aperiodic -- each $P$ carries exactly two periodic 
orbits. Moreover, the flow on $P$ is not standard on the boundary 
to account for the Dehn--twist in the surgery. Note that this method
also yields a flow on $M$ with exactly two periodic orbits.

The second result is that \emph{every compact three--manifold 
possesses a non-vanishing divergence--free $C^1$ vector 
field without periodic orbits}, \cite{kug}. In
particular, when applied to $S^3$, this theorem gives a 
volume--preserving
$C^1$ counterexample to the Seifert conjecture. The proof of the theorem
uses the previously mentioned result and a volume--preserving version 
of Schweitzer's plug to eliminate periodic orbits. To construct such a 
plug, G. Kuperberg applies, in a non-trivial way, stretching in the 
vertical direction to compensate for the area change resulting from
the Denjoy flow. (See \cite{kug} for more details.)

\begin{Remark}
\labell{rmk:Ratner}
There is a broad class of volume--preserving flows without periodic
orbits on homogeneous spaces of Lie groups. These flows arise as actions
of unipotent subgroups. 
For example, the horocycle flow (Section \ref{subsec:proofs}) and 
its Hamiltonian analogues from Example \ref{exam:hor-high} are
among such flows. Deep results on the closures of orbits for 
these flows and on their invariant measures are obtained by Ratner.
(See \cite{ratner} and references therein).
\end{Remark}

\section{Hamiltonian Vector Fields without Periodic Orbits}
\labell{sec:Ham}
The question of (non-)existence of periodic orbits for Hamiltonian
dynamical systems differs from that for general smooth dynamical
systems in at least one important way -- Hamiltonian systems tend
to have a lot of periodic orbits. For example, it is reasonable
to expect that a given Hamiltonian system with a proper smooth
Hamiltonian has a closed orbit on every regular energy level. 
Taken literally, this statement is not correct in general. (For example,
Zehnder, \cite{zehnder}, found a Hamiltonian flow on $\T^{2n}$,
$2n\geq 4$, with an irrational symplectic structure such that there
are no periodic orbits for a whole 
interval of energy values; see Example \ref{ex:Zehn} below.) However, 
periodic orbits are known 
to exist for almost all energy values for a broad class of symplectic
manifolds. For instance, as shown by
Hofer, Zehnder, and Struwe, \cite{ho-ze:per-sol,str},
almost all levels of a smooth proper
function on $\reals^{2n}$ carry at least one periodic orbit. 
A similar theorem holds for cotangent bundles, \cite{hv}. (The reader 
interested in a detailed discussion and further references 
should consult \cite{ho-ze:book}.) 

Furthermore, according to the $C^1$-closing lemma of Pugh and Robinson, 
\cite{pugh-rob}, a periodic orbit can be created from a recurrent 
trajectory by a $C^1$-small smooth perturbation of the original
dynamical system. (As Carlos Gutierrez pointed out to the author,
the perturbation can in fact be made $C^\infty$-smooth, 
\cite{gut:letter}.) As a consequence of the closing lemma, a 
$C^1$-generic system has the union of its periodic trajectories
dense in the set of its recurrent points. Both of these results hold
for Hamiltonian systems, \cite{pugh-rob}.

The situation becomes more subtle when the $C^1$-topology is replaced
by the $C^k$-topology with $k>1$. The $C^2$-closing lemma has not been
yet proved or disproved. (See, e.g., \cite{gut:example}, \cite{Car},
and \cite{AZ} for partial results, examples, and references.) 
In the Hamiltonian setting,
the systems whose periodic trajectories are dense may no longer be
generic if $k$ is roughly speaking greater than the dimension
of the manifold. For example, according to M. Herman, 
\cite{herman1,herman2}, Hamiltonian vector fields
with Hamiltonians $C^k$-close to Zehnder's Hamiltonian on 
$\T^{2n}$ do not have
periodic orbits for an interval of energy values, provided that
the symplectic form satisfies a certain Diophantine condition and $k>2n$.
These examples are, however, in some sense exceptional. In fact, the
theorems on the density of energy values for periodic orbits can be 
interpreted as that the Hamiltonian $C^k$-closing lemma holds in a very 
strong form for many symplectic manifolds.
It is still not known if the $C^k$-closing lemma with $k\geq 2$ fails 
in general for Hamiltonian flows when the symplectic form is exact near 
the energy level.

The examples of Hamiltonian flows without periodic orbits on
a compact energy level are scarce. In Section \ref{sec:list} we
attempt to list all known constructions of such flows.
In this section, we focus on the Hamiltonian Seifert conjecture
and on one particular method to construct Hamiltonian vector fields
without periodic orbits.

\begin{Remark}
There are a number of results concerning existence of 
periodic orbits on a fixed
energy level. Recall that according to Weinstein's conjecture, there
is a periodic orbit on an energy level of contact type, \cite{we:conj}.
This conjecture was proved for $\R^{2n}$ by Viterbo in \cite{vi:Theorem}
and then for many other symplectic manifolds (e.~g., for cotangent
bundles in \cite{hv}). The reader interested in details and more 
up-to-date references should consult \cite{ho}, \cite{ho-ze:book}, or 
\cite{vi:functors}.
Also note that in the context of contact topology and hydrodynamics the Seifert conjecture
is discussed in \cite{GE}.

\end{Remark}

\subsection{The Seifert Conjecture for Hamiltonian Vector Fields}
\labell{subsec:seif-ham}
The Seifert conjecture can be extended to Hamiltonian flows in a number
of ways. For example, one may ask if there is a proper smooth
function on a given symplectic manifold (e.g., $\reals^{2n})$, having
a regular level without periodic orbits. 
Recall that a \emph{characteristic} of a two-form $\eta$ of rank
$(2n-2)$ on a $(2n-1)$-dimensional manifold is an integral curve
of the field of directions formed by the null-spaces $\ker\eta$.
Thus, the question can be reformulated as whether or not in a given 
symplectic manifold there exists
a regular compact hypersurface without closed characteristics. One can
even ask whether the manifold admits a function with a sufficiently
big set of energy levels without periodic orbits. 

Let $i\colon M\hookrightarrow W$ be an embedded smooth compact 
hypersurface without boundary in a $2n$-dimensional symplectic manifold 
$(W,\sigma)$. 

\begin{Theorem}
\labell{Theorem:main1}
Assume that $2n\geq 6$ and that $i^*\sigma$ has only a finite number of 
closed characteristics. Then there exists a $C^\infty$-smooth embedding 
$i'\colon M\hookrightarrow W$, which is $C^0$-close and isotopic to $i$, 
such that ${i'}^*\sigma$ has no closed characteristics.
\end{Theorem}

An irrational ellipsoid $M$ in the standard symplectic vector space 
$\reals^{2n}=W$ corresponds to a collection of $n$ uncoupled harmonic
oscillators whose frequencies are linearly independent over 
$\rationals$, i.~e., $M$ is the unit level of a quadratic Hamiltonian 
with rationally independent frequencies. Thus, $M$ carries exactly $n$ 
periodic orbits. 
Applying Theorem \ref{Theorem:main1}, we obtain the following

\begin{Corollary}    
\labell{Corollary:sphere}
For $2n\geq 6$, there exists a $C^\infty$-embedding 
$S^{2n-1}\hookrightarrow \reals^{2n}$ such that the restriction of the 
standard symplectic form to $S^{2n-1}$ has no closed characteristics.
\end{Corollary}

\begin{Corollary}
\labell{Corollary:function}
For $2n\geq 6$, there exists a $C^\infty$-function 
$h\colon\reals^{2n}\to \reals$, $C^0$-close and isotopic
(with a compact support) to a positive definite quadratic form,
such that the Hamiltonian flow of $h$ has no closed trajectories on 
the level set $\{ h=1\}$.
\end{Corollary}

\begin{Remark}
These results, proved in \cite{gi:seifert} and, independently, 
by M. Herman \cite{herman-fax}, first
required the ambient space to be at least eight-dimensional, i.~e., 
$2n\geq 8$, in the $C^\infty$-case. A 
$C^{2+\epsilon}$-hypersurface
in $\reals^6$ was found by M. Herman \cite{herman-fax}. In 
\cite{gi:seifert97}, Theorem \ref{Theorem:main1} and its corollaries
were extended to $2n=6$.
\end{Remark}

\begin{Remark} Almost nothing is known on how big the
set of energy values $\CE$ of the levels without periodic orbits can be.
It is clear that as in Corollary \ref{Corollary:function}, there exists 
a function on $\reals^{2n}$ for which $\CE$ is infinite (but discrete).
It also seems plausible that $\CE$ can have limit points at critical
values of the Hamiltonian.
However, it is unknown whether or not $\CE$ can have a limit point 
that is a regular value. 

As is clear from the argument of \cite{gi:seifert,gi:seifert97}, there
exists a $C^0$-foliation of an open neighborhood of
an ellipsoid in $\R^{2n}$ such that every leaf is $C^\infty$-smooth and 
isotopic to the sphere and no leaf carries closed characteristics.
\end{Remark}

\begin{Remark} 
Theorem \ref{Theorem:main1} can be applied to any compact
hypersurface with a finite number of closed characteristics.
The only known examples of such hypersurfaces in $\R^{2n}$ 
that are not 
diffeomorphic to $S^{2n-1}$ are non-simply connected hypersurfaces 
constructed by Laudenbach, \cite{laud}.
\end{Remark}

An alternative way to formulate the Hamiltonian Seifert conjecture is 
to consider an odd--dimensional manifold with a maximally non-degenerate
closed two-form (rather than the restriction of
the symplectic form to a hypersurface).
A result similar to Theorem \ref{Theorem:main1}, involving only 
two-forms on $M$ but no symplectic embedding, is also correct. 

Let $\dim M=2n-1$
and let $\eta$ be a closed maximally non-degenerate (i.~e., of rank
$(2n-2)$) two-form on $M$.
Recall that two such forms $\eta$ and $\eta'$ are said to be \emph{homotopic} 
if there exists a family $\eta_\tau$, $\tau\in [0,1]$, 
of closed maximally non-degenerate forms connecting $\eta=\eta_0$ with 
$\eta'=\eta_1$ and such that all $\eta_\tau$ have the same cohomology 
class.  The following theorem
is proved in \cite{gi:seifert} and \cite{gi:seifert97}.

\begin{Theorem}
\labell{Theorem:main2}
Assume that $2n-1\geq 5$ and that $\eta$ has a finite number of 
closed characteristics. Then there exists a closed maximally 
non-degenerate 2-form $\eta'$ on $M$ which is homotopic to $\eta$ and has 
no closed characteristics.
\end{Theorem}

\begin{Remark}
\labell{rmk:hms}
In fact, Theorem \ref{Theorem:main2} is a corollary of Theorem 
\ref{Theorem:main1}. To derive Theorem \ref{Theorem:main2} from Theorem 
\ref{Theorem:main1}, consider a pair $(M^{2n-1},\eta)$, where
$\eta$ is closed and maximally non-degenerate. Then, as follows from
Gotay's coisotropic embedding theorem, \cite{gotay}, there exists a 
symplectic manifold $(W^{2n},\sigma)$ and a proper embedding 
$M\hookrightarrow W$ such that $\sigma|_M=\eta$.

More explicitly, set $W=M\times (-1,1)$ and let $t$ be the coordinate 
on $(-1,1)$.
To construct $\sigma$, fix a 1-form $\alpha$ on $M$ such that 
$\ker \alpha$ is everywhere transversal to the characteristics of $\eta$.
Then $\sigma=\eta+d(\epsilon t\alpha)$ is symplectic on $W$, provided
that $\epsilon>0$ is small enough. (The author is grateful to Ana Cannas 
da Silva for this remark.)

Technically, however, it is more convenient to prove
Theorem \ref{Theorem:main2} first and then to modify its proof
to obtain Theorem \ref{Theorem:main1}.
\end{Remark}

\begin{Remark}
Theorem \ref{Theorem:main2} extends to the real analytic case: one can
make the form $\eta'$ real analytic, provided that $Q$ and $\eta$ are 
real analytic. The argument is the same as that used in the construction 
of a real analytic version of Wilson's, \cite{wilson}, or Kuperberg's, 
\cite{kuk}, flow. (See \cite{ghys} and \cite{kugk}.)
\end{Remark}

\subsection{The Hamiltonian Seifert Theorem} Before we turn to the
outline of the proofs of Theorems \ref{Theorem:main1} and 
\ref{Theorem:main2}, let us state a result which can be viewed as a
Hamiltonian version of Theorem \ref{Theorem:fuller}.

Recall that a periodic orbit of a vector field or an integral curve
of a one-dimensional foliation is called non-degenerate if the 
linearization of its Poincar\'{e} return map does not have unit as
an eigenvalue. Let $(B, \omega)$ be a compact symplectic manifold and let
$\pi\colon M\to B$ be a principle $S^1$-bundle. For a one-form $\lambda$
such that $d\lambda$ is $C^0$ sufficiently small, the closed two-form 
$\eta=\pi^*\omega+d\lambda$ on $M$ is maximally non-degenerate (and
homotopic to $\pi^*\omega$).

\begin{Theorem}[\cite{gi:MathZ}]
\labell{thm:Ham-seifert}
If $d\lambda$ is $C^0$-small enough, the number of 
closed characteristics
of $\eta$ is strictly greater than the cup-length of $B$. If, 
in addition, all closed characteristics are non-degenerate, then the number
of closed characteristics is greater than or equal to the sum of 
Betti numbers of $B$.
\end{Theorem}

Theorem \ref{thm:Ham-seifert} generalizes some of the results (for
equal eigenvalues) of \cite{weinstein-1973,Moser-1976,Bottkol}.
A simple geometrical proof for the case where $B$ is a surface and 
$\lambda$ is $C^2$-small can be found in \cite{gi:FA}. 
Theorem \ref{thm:Ham-seifert} is closely related to the Arnold 
conjecture and to the problem of existence of periodic orbits of
a charge in a magnetic field; \cite{gi:Cambr}. 
Note also that if $d\lambda$ is not $C^0$-small, but
$\eta$ is still non-degenerate and homotopic to $\pi^*\omega$, the form 
$\eta$ may have no closed characteristics at all. The example is the
horocycle flow described in the next section. (See \cite{gi:Cambr} 
and \cite{gi:MathZ} for more details.)

\subsection{Proofs of Theorems \ref{Theorem:main1} and \ref{Theorem:main2}
and Symplectic Embeddings.}
\labell{subsec:proofs}
There are two essentially different methods to prove these theorems.
Both approaches follow the same general line as Wilson's argument
(i.~e., the proof of Theorem \ref{thm:wilson}), modified to make it work 
in the symplectic category, and vary only in the dynamics of the plug. 

The plug introduced by M. Herman, \cite{herman-fax}, is the 
symplectization
of Wilson's plug when $2n\geq 8$ or of the Schweitzer--Harrison plug
when $2n=6$ (hence the $C^{2+\eps}$-smoothness constraint) modified at 
infinity. In other words, these plugs are obtained by taking the 
induced flows on the cotangent bundles to the non-Hamiltonian plugs
and altering these flows away from a neighborhood of the zero section while
keeping all the properties of the plug.

The plugs from \cite{gi:seifert,gi:seifert97}
are built similarly to Wilson's plug but the entire construction
is Hamiltonian. Let us give some more details on this method focusing
specifically on the proof of Theorem \ref{Theorem:main2}. (The passage 
to Theorem \ref{Theorem:main1} is simply by showing that the same 
argument can be carried out for hypersurfaces. See \cite{gi:seifert}.) 

In the construction of the plug $P$ below we identify a neighborhood 
in $M$ containing $P$ with a small ball in $\reals^{2n-1}$ 
equipped with the two-form $\sigma$ induced by the canonical inclusion 
$\reals^{2n-1}\subset\reals^{2n}$.

The key element that makes Wilson's construction work in the 
Hamiltonian category is that the core of the plug (i.~e., the flow 
replacing the irrational flow on the torus)
is also chosen to be Hamiltonian. In other words, to have the aperiodicity 
condition satisfied, we need to find a Hamiltonian flow on a symplectic 
manifold having no periodic orbits on some compact energy level. 
Furthermore, to have the embedding condition satisfied, this energy
level with the induced two-form has to be embeddable into 
$\reals^{2n-1}$ to be made into a part of the plug. 
(In particular, the symplectic form must be exact on a neighborhood 
of the energy level.) The Hamiltonian flow used 
\cite{gi:seifert,gi:seifert97} is the horocycle flow, which we will 
now describe.

Let $\Sigma$ be a closed surface with a hyperbolic metric, i.~e., a
metric with constant negative curvature $K=-1$. Denote by
$\lambda$ the canonical Liouville one-form ``$p\,dq$'' on $T^*\Sigma$
and by $\Omega$ the pull-back to $T^*\Sigma$ of the area form on $\Sigma$.
The form $\omega=d\lambda+\Omega$ is symplectic on $T^*\Sigma$ and
exact on the complement of the zero section. The Hamiltonian flow 
of the standard metric Hamiltonian on $T^*\Sigma$ with respect to
the twisted symplectic form $\omega$ is known to have no periodic
orbits on the unit energy level $ST^*\Sigma$. Indeed,
the restriction $\varphi^t$ of this flow to $ST^*\Sigma$ is the standard 
horocycle
flow. The fact that $\varphi^t$ has no periodic orbits is a consequence of
the classical result of Hedlund, \cite{hedlund}, that the horocycle flow
is minimal, i.~e., all orbits are dense. We will return to the horocycle
flow in Section \ref{sec:list}.

To construct the plug, we then need to prove that there is a 
``symplectic'' embedding of a neighborhood $U$ of $ST^*\Sigma$ into 
$\reals^{2n-1}$, i.~e., an embedding $j\colon U\hookrightarrow\R^{2n-1}$ 
such that $j^*\sigma=\omega$. When $2n-1\geq 7$, this readily follows 
from a general symplectic embedding theorem due to Gromov, 
\cite{gr:icm,gr:book}. 

When $2n-1=5$ an additional argument is required. Consider the forms
$\omega_t=d\lambda+t\Omega$ on $T^*\Sigma$. These forms are
symplectic for all $t$. The form $\omega_0$ is the standard symplectic 
form on $T^*\Sigma$ and $\omega_1$ is the twisted form $\omega$.
Note that if $U$ is small enough, there exists a ``symplectic''
embedding $j_0\colon (U,\omega_0)\hook (\reals^5, \sigma)$. This is an
easy consequence of the fact that
there exists a Lagrangian immersion $\Sigma\to\reals^4$, 
\cite{gr:icm,gr:book}. Note also that $U$ can be replaced
by a closed neighborhood of $ST^*\Sigma$. Then the existence of a 
symplectic
embedding for $2n-1=5$ is a particular case of the following 
result improving the dimensional constraints from \cite{gr:icm,gr:book} 
by one.

Let $U$ and $V$ be manifolds of
equal dimensions. The manifold $U$ is assumed to be compact, perhaps with
boundary, while $V$ may be open but must be a manifold without 
boundary. Let $\sigma$ be a symplectic form on $V$. Abusing notation,
also denote by $\sigma$ the pull-back of $\sigma$ to $V\times\reals$ 
under the natural projection $V\times\reals\to V$.

\begin{Theorem}[\cite{gi:seifert97}]
\labell{Theorem:embed}
Let $\omega_t$, $t\in [0,1]$, be a family of symplectic forms
on $U$ in a fixed cohomology class: $[\omega_t]=\const$. Assume
that there is an embedding 
$j_0\colon U\hookrightarrow V\times\reals$ such
that $j_0^*\sigma=\omega_0$. Then there exists an embedding 
$j\colon U\hookrightarrow V\times\reals$, isotopic to $j_0$, with
$j^*\sigma=\omega_1$.
\end{Theorem}

\begin{Remark}
\labell{Remark:moser1}
Since $\sigma|_V$ is symplectic, the composition of $j$ with
the projection to $V$ is necessarily an immersion.
When $\p U=\emptyset$, Theorem \ref{Theorem:embed} follows immediately
from Moser's theorem \cite{moser}. 
\end{Remark}

The rest of the proof (i.~e., the construction of the form $\eta'$
on the complement of the core) proceeds according to the same scheme 
as the proof of 
Theorem \ref{thm:wilson} with more or less straightforward modifications;
see \cite{gi:seifert} for details.
This completes the outline of the proofs of Theorems \ref{Theorem:main1} 
and \ref{Theorem:main2}.

\section{The List of Hamiltonian Flows Without Periodic Orbits.} 
\labell{sec:list}

In this section we list all known to the author examples of smooth 
Hamiltonian systems on symplectic manifolds $(V,\omega)$ having
no periodic orbits on a compact regular level $M$. This list is similar to
the one given in \cite{gi:seifert97} with the exception of Example
\ref{exam:hor-high}. 

The list is divided into two parts according to whether $\omega$
is exact near $M$ or not. The reason for this division is that 
the qualitative behavior of the flows for which $[\omega|_M]\neq 0$ can be
expected to differ in an essential way from the behavior of the flows 
with $[\omega|_M]=0$. 
For instance, Zehnder's flow (Example \ref{ex:Zehn}) appears to be more 
robust than it seems to be possible for
``exact'' flows.

\bigskip
\noindent {\bf Case 1:}{\em~ The form $\omega$ is not exact.}

\begin{Example}[Zehnder's torus, \cite{zehnder}]
\labell{ex:Zehn}
Let $2n\geq 4$. Consider the torus $V=\TT^{2n}$ with an
irrational translation-invariant symplectic structure $\omega$.
Choose a Hamiltonian $H$ on $V$ so that every level $\{ H=c\}$ with
$c\in (0.5, 1.5)$ is the union 
of two standard embedded tori $\TT^{2n-1}\subset \TT^{2n}$. Since $\omega$
is irrational, the characteristics of $\omega|_{\TT^{2n-1}}$ form an
quasiperiodic flow on $\TT^{2n-1}$. Thus, none of the levels
$\{H=c\}$ with $c\in (0.5, 1.5)$ carries a periodic orbit. 
Note that $\omega$ is not exact on any of these energy levels.
As  we have already pointed out,
the flow in question exhibits remarkable stability properties according
to a result of M. Herman, \cite{herman1,herman2}. 
\end{Example}

\begin{Example}[The Hamiltonian horocycle flow in dimension $2n$] 
\labell{exam:hor-high}
Let $\C H^n$ $\,$ be the complex hyperbolic space equipped with its
standard K\"{a}hler metric (see, e.g., Section XI.10 of \cite{kob-nom}).
Pick a discrete subgroup $\Gamma$ in the group $\mbox{\rm SU}(1,n)$ 
of Hermitian isometries of $\C H^n$ such that $B=\Gamma\backslash\C H^n$ 
is smooth and compact. (To see that $\Gamma$ exists recall that 
according to a theorem of Borel, \cite{borel}, there is a discrete 
subgroup $\Gamma$ such that $\Gamma\backslash\C H^n$ is compact. Then 
by Selberg's lemma, \cite{selberg}, $\Gamma$ can be chosen
so that $B$ is smooth.) 

Let $H$ be the standard metric Hamiltonian and let $d\lambda$ be
the standard symplectic structure on $T^*B$. Denote by $\Omega$ 
the pull-back to $T^*B$ of the K\"{a}hler symplectic form on $B$. 
The Hamiltonian flow of $H$ with respect to
the twisted symplectic structure $\omega=\Omega+d\lambda$ has no periodic
orbits on the level $M=\{H=1\}$. This follows from the fact that this 
flow is
generated by the (right) action of a unipotent one-parameter subgroup
of $\mbox{\rm SU}(1,n)$ and that $\Gamma$ contains no unipotent elements, 
\cite{KM}; see also \cite{rag}, \cite{ratner}, and references therein.

When $n=1$ this construction gives exactly the horocycle flow on the
unit (co)tangent bundle to a surface $B=\Sigma$. If $n>1$, the form
$\omega$ is not exact in any neighborhood of $M$.

Note also that the Hamiltonian horocycle flows arise in the 
description of the motion of a charge in a magnetic field on the
configuration space $B$; see \cite{gi:Cambr,gi:MathZ,ely:paper}.
These are the only known ``magnetic'' flows without 
periodic orbits on an energy level.

\end{Example}

\noindent {\bf Case 2:}{\em~ The form $\omega$ is exact near the energy level.}
In this case we should distinguish whether $\dim V=4$ or $\dim V \geq 6$.

When $\dim V=4$, the only known smooth example is the horocycle 
flow described as a Hamiltonian system in Section 
\ref{subsec:proofs} and in Example \ref{exam:hor-high} with $n=1$ and 
$B=\Sigma$. To slightly generalize this
construction, note that a neighborhood of $ST^*\Sigma$ can be identified 
with
$
V=\Gamma\backslash \mbox{\rm SU}(1,1)\times (1-\epsilon,1+\epsilon)
$,
where 
$\Gamma=\pi_1(\Sigma)$, so that $H$ becomes the projection to the second 
component. Then, instead of $\Gamma=\pi_1(\Sigma)$ 
we can take any 
discrete subgroup such that the quotient 
$\Gamma\backslash\mbox{\rm SU}(1,1)$ 
is compact and smooth. As follows from Remark \ref{rmk:hms}, the flow on 
$\Gamma\backslash\mbox{\rm SU}(1,1)$ 
generated by the action of a unipotent subgroup of 
$\mbox{\rm SU}(1,1)$ is the Hamiltonian flow of $H$ on $\{ H=1 \}$ 
with respect to some symplectic form. 

Finally note that by Remark \ref{rmk:hms} the flows on three--manifolds 
constructed by G. Kuperberg and described above in Section 
\ref{subsec:vol} can also be thought of as Hamiltonian flows on some
symplectic manifolds. This gives a class of smooth Hamiltonian flows
on symplectic manifolds with a finite number of closed orbits on a given 
energy level or $C^1$-flows without periodic orbits. It is not known
if such a flow on $S^3$ can be obtained by a $C^2$-embedding of
$S^3$ into $\R^4$.

When $\dim V\geq 6$, the only known examples are essentially those
described in Section \ref{subsec:seif-ham} or those obtained by
iterations of their constructions. In other words, beginning with a
flow with a finite number of periodic orbits one can eliminate these
orbits by using the plugs from \cite{gi:seifert,gi:seifert97} or 
\cite{herman-fax}. The resulting manifold can also be used as the 
core of a plug in the same way as the horocycle flow. As in
Remark \ref{rmk:proliferation}, these new 
plugs can in turn be employed to construct Hamiltonian flows without 
periodic orbits, etc.

Since the starting point of this method is a flow with
a finite number of periodic orbits, it makes sense to ask how many
such flows are known in addition to those on irrational ellipsoids.
For example, we have already mentioned non-simply connected 
hypersurfaces in $\R^{2n}$ found by Laudenbach \cite{laud}. 
To produce more examples, we can apply the same method as described
in the previous paragraph. Namely,
if a Hamiltonian flow with a finite number of periodic orbits can be
used as the core of the plug, the resulting ``plug'' will also
have only a finite number of periodic orbits inside of it. Thus by
inserting it into another flow with a finite number of periodic orbits,
we can create yet one more flow with the same property. For example, by 
taking $S^1\subset \R^2$ as the core, Cieliebak, \cite{ciel}, constructed 
embeddings $S^{2n-1}\subset\reals^{2n}$, $2n\geq 4$, such that the pairs
of closed characteristics are linked and knotted in an essentially 
arbitrary way. (Each plug results into two ``parallel'' closed orbits,
but there are no constraints on knotting of these pairs or their linking.)

\end{document}